%% file: dipole.tex
\documentclass{article}
\usepackage{dipole}

\begin{document}
\title{Bipolar comparison}
\author{Nina Lebedeva, Anton Petrunin and Vladimir Zolotov}

\nofootnote{The first author was partially supported by RFBR grant 17-01-00128, 
the second author was partially supported by NSF grant DMS 1309340,
the third author was supported in part by DFG grant SPP 2026 and RFBR grant 17-01-00128.
Essential part of this work was done during the intense activity period ``Metric Measure Spaces and Ricci Curvature'' September 4--29, 2017 in Bonn.}

\newcommand{\Addresses}{{\bigskip\footnotesize

  Nina Lebedeva, \par\nopagebreak
  \textsc{Steklov Institute,
27 Fontanka, St. Petersburg, 
191023, Russia.}
  \par\nopagebreak
  \textit{Email}: \texttt{lebed@pdmi.ras.ru}

\medskip

  Anton Petrunin, \par\nopagebreak\textsc{Math. Dept. PSU, University Park, PA 16802, USA.}
  \par\nopagebreak
  \textit{Email}: \texttt{petrunin@math.psu.edu}
  
\medskip

Vladimir Zolotov, \par\nopagebreak\textsc{Steklov Institute,
27 Fontanka, St. Petersburg, 
191023, Russia.}
  \par\nopagebreak
  \textit{Email}: \texttt{paranuel@mail.ru}
}}

\date{}

\maketitle

\begin{abstract}
We define a new type of metric comparison similar to the comparison of Alexandrov.
We show that it has strong connections to continuity of optimal transport between regular measures on a Riemannian manifold, in particular to the so called MTW condition introduced by Xi-Nan Ma, Neil Trudinger and Xu-Jia Wang.

\end{abstract}

\input{defs.tex}

\input{kirszbraun.tex}
\input{6-dipole.tex}
\input{convexity.tex}

\input{mtw.tex}

\input{polypole.tex}
\input{end.tex}

\Addresses
\end{document}

%% file: defs.tex
\section{Introduction}\label{sec:intro}

We will denote by $|a-b|_X$ the distance between points $a$ and $b$ in the metric space $X$.

\parbf{Tree comparison.}
Fix a tree $T$ with $n$ vertexes.

Let $(a_1,\dots, a_n)$ be a point array in a metric space $X$ labeled by the vertexes of $T$.
We say that $(a_1,\dots, a_n)$  satisfies the \emph{$T$-tree comparison} if there is a point array $(\~a_1,\dots, \~a_n)$ in the Hilbert space $\HH$ such that 
\[|\~a_i-\~a_j|_\HH\ge|a_i-a_j|_X\]
for any $i$ and $j$ and equality holds if $a_i$ and $a_j$ are adjacent in $T$.

We say that a metric space $X$ satisfies the \emph{$T$-tree comparison} if 
every $n$-points array in $X$ satisfies the $T$-tree comparison.

Instead of the Hilbert space $\HH$, we may use an infinite dimensional sphere or an infinite dimensional hyperbolic space.
In this case it defines \emph{spherical} and \emph{hyperbolic} tree comparisons.

\begin{wrapfigure}{r}{27 mm}
\vskip-8mm
\begin{tikzpicture}[scale=1,
  thick,main node/.style={circle,draw,font=\sffamily\bfseries,minimum size=3mm}]
  \node[main node] (1) at (-1/2,-5/6) {$x$};
  \node[main node] (2) at (0,0){$p$};
  \node[main node] (3) at (-1/2,5/6){$y$};
  \node[main node] (4) at (3/2,5/6) {$v$};
  \node[main node] (5) at (1,0) {$q$};
  \node[main node] (6) at (3/2,-5/6) {$w$};

  \path[every node/.style={font=\sffamily\small}]
   (1) edge node[above]{}(2)
   (2) edge node[above]{}(3)
   (2) edge node[above]{}(5)
   (4) edge node[above]{}(5)
   (5) edge node[above]{}(6);
\end{tikzpicture}
\end{wrapfigure}

\parbf{Encoding of trees.}
To encode the labeled tree on the diagram, we will use notation $p/xy(q/vw)$.
It means that we choose $p$ as the root; 
$p$ has two children leaves $x$, $y$ and one child $q$ with two children leaves $v$ and $w$.
Taking another root for the same tree, we get different encodings, for example $q/vw(p/xy)$ or $x/(p/y(q/vw))$.

If we do not need the labeling of vertexes,
it is sufficient to write the number of leaves in the brackets;
this way we can write 2(2) instead of $p/xy(q/vw)$ since the root ($p$) has 2 leaves ($x$ and $y$) and yet another child ($q$) that has 2 leaves ($v$ and $w$).  
The same tree can be encoded as (1(2)) meaning that the root $x$ has no leaves, 
$p$ has 1 leaf $y$ and one child $q$ with 2 leaves $v$ and $w$.
Note that every vertex that is not the root and not a leaf corresponds to a pair of brackets in this notation.

Using the described notation, we could say that a metric space \emph{satisfies the 2(2)-tree comparison},  meaning that it satisfies the tree comparison on the diagram.
We could also say \emph{``applying the tree comparison for $p/xy(q/vw)$ ...''} meaning that we apply the comparison for these 6 points in a metric space labeled as on the diagram.

\parbf{Monopolar trees.}
A vertex of a tree of degree at least two will be called a \emph{pole}.

\begin{wrapfigure}{r}{21 mm}
\begin{tikzpicture}[scale=1,
  thick,main node/.style={circle,draw,font=\sffamily\bfseries,minimum size=3mm}]

  \node[main node] (1) at (5/6,1) {};
  \node[main node] (2) at (0,3/2){};
  \node[main node] (3) at (10/6,3/2){};
  \node[main node] (4) at (5/6,0) {};

  \path[every node/.style={font=\sffamily\small}]
   (1) edge node[above]{}(2)
   (1) edge node[above]{}(3)
   (1) edge node[above]{}(4);
\end{tikzpicture}
\end{wrapfigure}

Recall that \emph{Alexandrov space} with nonnegative curvature is defined as a complete length space with \emph{nonnegative curvature in the sense of Alexandrov};
the latter is equivalent to the 3-tree comparison; that is, the comparison for the tripod-tree on the diagram. 

Using the introduced notation, the theorem on ($n$+1)-point comparison in \cite{AKP} can be restated the following way: \emph{If a complete length-metric space satisfies $3$-tree comparison, then it also satisfies $n$-tree comparison for every positive integer~$n$; in other words it satisfies all monopolar tree comparisons.}

\parbf{Bipolar comparison.}
The following theorem is proved in sections~\ref{sec:pivotal} and \ref{6-dipole};
it describes the comparisons for the bipolar trees 3(1) and 2(2) shown on the diagram.

\begin{center}
\begin{tikzpicture}[scale=1,
  thick,main node/.style={circle,draw,font=\sffamily\bfseries,minimum size=3mm}]

  \node[main node] (1) at (1,2) {};
  \node[main node] (2) at (1,0){};
  \node[main node] (3) at (1,1){};
  \node[main node] (4) at (0,1) {};
  \node[main node] (5) at (3,1) {};
  \node[main node] (6) at (2,1) {};

  \path[every node/.style={font=\sffamily\small}]
   (1) edge node[above]{}(3)
   (2) edge node[above]{}(3)
   (3) edge node[above]{}(6)
   (4) edge node[above]{}(3)
   (5) edge node[above]{}(6);
\end{tikzpicture}
\hskip30mm
\begin{tikzpicture}[scale=1,
  thick,main node/.style={circle,draw,font=\sffamily\bfseries,minimum size=3mm}]

  \node[main node] (1) at (0,0) {};
  \node[main node] (2) at (1/2,5/6){};
  \node[main node] (3) at (0,10/6){};
  \node[main node] (4) at (2,0) {};
  \node[main node] (5) at (3/2,5/6) {};
  \node[main node] (6) at (2,10/6) {};

  \path[every node/.style={font=\sffamily\small}]
   (1) edge node[above]{}(2)
   (2) edge node[above]{}(3)
   (2) edge node[above]{}(5)
   (4) edge node[above]{}(5)
   (5) edge node[above]{}(6);
\end{tikzpicture}
\end{center}

\begin{thm}{Theorem}\label{thm:3(1)+2(2)}
A complete Riemannian manifold has nonnegative sectional curvature if and only if it satisfies 3(1)-tree (or 2(2)-tree) comparison.
\end{thm}

{

\begin{wrapfigure}{r}{32 mm}
\vskip-2mm
\begin{tikzpicture}[scale=1,
  thick,main node/.style={circle,draw,font=\sffamily\bfseries,minimum size=3mm}]

  \node[main node] (0) at (.3,.9){};
   \node[main node] (1) at (-.8,-.6){};
  \node[main node] (2) at (.3,-.9){};
  \node[main node] (3) at (0,0){};
  \node[main node] (4) at (-.8,.6){};
  \node[main node] (5) at (2,0){};
  \node[main node] (6) at (1,0){};

  \path[every node/.style={font=\sffamily\small}]
     (0) edge node[above]{}(3)
   (1) edge node[above]{}(3)
   (2) edge node[above]{}(3)
   (3) edge node[above]{}(6)
   (4) edge node[above]{}(3)
   (5) edge node[above]{}(6);
\end{tikzpicture}
\end{wrapfigure}

The proof is given in Section~\ref{6-dipole}.

\bigskip

The following theorem gives a description of Riemannian manifolds satisfying 4(1)-tree comparison (the 4(1)-tree is shown on the diagram).
Let us first define CTIL Riemannian manifolds (CTIL stands for \emph{convexity of tangent injectivity locus}).

}

Let $M$ be a Riemannian manifold and $p\in M$.
The subset of tangent vectors $v\in\T_p$ such that there is a minimizing geodesic $[p\,q]$ in the direction of $v$ with length $|v|$ will be denoted as $\overline{\TIL}_p$.
The interior of $\overline{\TIL}_p$ is denoted by $\TIL_p$; it is called \emph{tangent injectivity locus} at $p$.
If at $\TIL_p$ is convex for any $p\in M$, then $M$ is called CTIL.

For a function $f$ defined on an open convex set of Euclidean space we write 
$f''\le \lambda$ if for any unit-speed geodesic $\gamma$ the function
\[t\mapsto f\circ\gamma(t)-\tfrac\lambda2\cdot t^2\]
is a concave real-to-real function.

\begin{thm}{Main theorem}\label{thm:convexity}
If a complete Riemannian manifold $M$ satisfies $4(1)$-tree comparison, then it is CTIL and the following two conditions hold:
\begin{enumerate}[(i)]
\item\label{thm:convexity:convexity} For any $p,q\in M$, we have $f''\le 1$, where $f$ is the function $f\: \TIL_p\to \RR$ defined by
\[f(v)=\tfrac12\cdot\dist_q^2\circ\exp_p(v).\] 
\item\label{thm:convexity:MTW} For any point $p\in M$ and any three tangent vectors 
$W\in \TIL_p$, $X,Y\in \T_p$, we have
\[\frac{\partial^4}{\partial^2s\,\partial^2t}\left|\exp_p(s\cdot X)-\exp_p(W+t\cdot Y)\right|_M^2\le 0\eqlbl{eq:MTW}\]
at $t=s=0$.
\end{enumerate}

Moreover, if one of the conditions (\ref{thm:convexity:convexity}) or (\ref{thm:convexity:MTW}) holds in a CTIL Riemannian manifold $M$, then $M$ satisfies all bipolar comparisons; that is, the $m(n)$-tree comparison holds for any $m$ and $n$.
\end{thm}

The part of theorem related to the condition (\textit{\ref{thm:convexity:convexity}}) is proved in Section~\ref{convexity}.
The equivalence of the conditions (\textit{\ref{thm:convexity:convexity}}) and (\textit{\ref{thm:convexity:MTW}}) for CTIL manifolds is proved in Section~\ref{MTW+}.
These proofs and the proof of Theorem~\ref{thm:3(1)+2(2)} are essentially independent; a possible generalization to length-metric spaces is discussed in the final remarks.

Note that theorems \ref{thm:3(1)+2(2)} and \ref{thm:convexity} do not provide descriptions of manifolds only for the following two bipolar trees: 3(2) and 3(3).

A Riemannian manifold $M$ is called MTW if inequality \ref{eq:MTW} holds for any $p\in M$ and all triples of vectors $W\in \TIL_p$, $X,Y\in \T_p$ with an additional assumption $X\perp Y$.
The condition (\textit{\ref{thm:convexity:MTW}}) in the theorem is stronger since it does not require orthogonality; this condition is named MTW$^{\not\perp}$ in \cite{FRV-Nec+Suf}.

The MTW condition is named for Xi-Nan Ma, Neil Trudinger and Xu-Jia Wang who introduced it in \cite{MTW};
an important step in the understanding MTW was made by Grégoire Loeper in~\cite{loeper}.
Below we describe a connection between MTW and the so called \emph{transport continuity property}, briefly TCP.

A compact Riemannian manifold $M$ is called TCP 
if for any two regular measures with density functions bounded away from zero and infinity the generalized solution of Monge--Amp\`{e}re equation provided by optimal transport 
is a genuine (continuous) solution.

MTW turns out to be a necessary condition for TCP.
In \cite{FRV-Nec+Suf}, 
Alessio Figalli, 
Ludovic Rifford 
and C\'edric Villani showed that
a strict version of CTIL and MTW provides a sufficient condition for TCP.
From the proof of Theorem~\ref{thm:convexity} it is evident that spherical 4(1)-tree comparison implies the strict version of MTW; in particular, we get the following:

\begin{thm}{Corollary}
Any Riemannian manifold satisfying spherical 4(1)-tree comparison is TCP.
\end{thm}

Note that the identity $s''\equiv1$ holds for the function $s\:v\mapsto\tfrac12\cdot|v|^2$ defined on $\TIL_p$.
Therefore the condition $f''\le 1$ in (\textit{\ref{thm:convexity:convexity}}) is equivalent to 
concavity of the function $h\zz=f-s$.
The following proposition gives a description of the MTW condition via convexity of superlevel set of $h$.
Essentially, this proposition is just a reformulation of a synthetic description of MTW given by C\'edric Villani \cite[Proposition 2.6]{MTW+CTIL}.
(Our formulation can be proved by recursive application of Villani's description and Villani's description follows immideately from our fromulation.) 
The proof can be also obtained by a straightforward modification of the proof of Theorem~\ref{thm:convexity}.

\begin{thm}{Proposition}\label{prop:convexity}
A complete CTIL Riemannian manifold $M$ is MTW if and only if 
for any $p,q\in M$, the function $h\: \TIL_p\to \RR$ defined by
\[h(v)=\tfrac12\cdot(\dist_q^2\circ\exp_p(v)-|v|^2)\] 
has convex superlevel sets; that is, the set
\[\set{v\in\TIL_p}{h(v)> C}\]
is convex for any $C\in\RR$.
\end{thm}

\parbf{All tree comparisons.}
Recall that a map $f\:W\to X$ between metric spaces is called \emph{submetry} if for any $w\in W$ and $r\ge 0$, we have 
\[f[B(w,r)_W]=B(f(w),r)_X,\]
where $B(w,r)_W$ denotes the open ball with center $w$ and radius $r$ in the space~$W$.
In other words submetry is a map that is 1-Lipschitz and 1-co-Lipschitz at the same time.
Note that by the definition, any submetry is onto.

\begin{thm}{Exercise}\label{ex:quotient}
Let $W\to X$ be a submetry and $T$ be a tree.
Assume $W$ satisfies the $T$-tree comparison.
Show that the same holds for $X$.
\end{thm}

Note that the Hilbert space satisfies all tree comparisons; it follows directly from the definition.
According to the exercise, the same holds for the target spaces of submetries defined of Hilbert space or its subsets.
The following theorem gives a converse of the latter statement.

\begin{thm}{Theorem}\label{thm:hilbert-quotient}
A separable metric space $X$ satisfies all tree comparisons if and only if
$X$ is isometric to a target space of submetry defined on a subset  of a separable Hilbert space.
\end{thm} 

The following proposition provides a source of examples of spaces satisfying all tree comparisons.
For example, since $\SS^n=\SO(n)/\SO(n-1)$, any round sphere has this property.

\begin{thm}{Proposition}\label{prop:group}
Suppose $G$ is a compact Lie group with bi-invariant metric, so the action $G\times G\acts G$ defined by $(h_1,h_2)\cdot g=h_1\cdot g\cdot  h_2^{-1}$ is isometric. 
Then for any closed subgroup $H<G\times G$, the bi-quotient space $G/\!\!/H$ satisfies all tree comparisons.
\end{thm}

The theorem and proposition are proved in Section~\ref{sec:all-tree}.

The proposition and Theorem~\ref{thm:convexity} imply that every bi-quotient space $G/\!\!/H$ is CTIL and MTW.
These examples seem to be new; a related source of examples is found by Young-Heon Kim and Robert McCann \cite{kim-mccann}.

%% file: kirszbraun.tex
\section{Kirszbraun's rigidity}

In the proof we will use the rigidity case of the generalized Kirszbraun theorem proved by Urs Lang and Viktor Schroeder in \cite{LS}, see also \cite{AKP}.

\begin{thm}{Kirszbraun rigidity theorem}\label{thm:kirszbraun-rigid}
Let $M$ be a complete Riemannian manifold with nonnegative sectional curvature.

Assume that for two point arrays $p,x_1,\dots,x_n\in M$ and $\~q, \~x_1,\dots,\~x_n\in \HH$ we have that 
\[|\~q-\~x_i|\ge |p-x_i|\]
for any $i$,
\[|\~x_i-\~x_j|\le |x_i-x_i|\]
for any pair $(i,j)$
and $\~q$ lies in the interior of the convex hull $\~K$ of $\~x_1,\dots,\~x_n$.

Then equalities hold in all the inequalities above.
Moreover there is an distance preserving map $f\:\~K\to M$ such that $f(\~x_i)=x_i$ and $f(\~q)=p$. 
\end{thm}

We reduce the theorem to the case $M=\RR^m$ which is left as an exercise.

\parit{Proof.}
By the generalized Kirszbraun theorem, there is a short map\footnote{Here and further \emph{short map} is a \emph{distance nonexpanding map}.} $f\:M\to \HH$
such that $f(x_i)=\~x_i$.
Set  $\~p=f(p)$.
By assumptions
\[|\~q-\~x_i|\ge |\~p-\~x_i|.\]

Since $\~q$ lies in the interior of $K$, we have that $\~q=\~p$ and the equality 
\[|\~q-\~x_i|= |p-x_i|.\]
holds for each $i$.

Set $v_i=\log_px_i$; that is, $t\mapsto \exp_p(t\cdot v_i)$ for $t\in[0,1]$ is a minimizing geodesic from $p$ to $x_i$ or, equivalently, $|v_i|=|p-x_i|$ and $\exp_pv_i=x_i$.
Recall that the gradient exponent $\gexp_p\:\T_p\to M$ is defined for any compete Riemannian manifold;
it is is a short map if $M$ has nonnegative curvature 
and $\gexp_p\:v_i\mapsto x_i$ for each $i$ (see \cite{AKP}).

The composition $f\circ \gexp_p\:\T_p\to \HH$ is short
and by classical Kirszbraun rigidity it has to be distance preserving on the convex hull $\~K'$ of $v_i$.
Hence $\~K'$ is isometric to $\~K$ and the restriction $g|_{\~K'}$ is distance preserving. 
Hence the result.
\qeds

%% file: 6-dipole.tex
\pagebreak 
\section{Pivotal trees}\label{sec:pivotal} 

\begin{wrapfigure}{r}{27 mm}
\vskip-0mm
\begin{tikzpicture}[scale=1,
  thick,main node/.style={circle,draw,font=\sffamily\bfseries,minimum size=3mm}]
  \node[main node] (1) at (-1/2,-5/6) {$x$};
  \node[main node] (2) at (0,0){$p$};
  \node[main node] (3) at (-1/2,5/6){$y$};
  \node[main node] (4) at (3/2,5/6) {$v$};
  \node[main node] (5) at (1,0) {$q$};
  \node[main node] (6) at (3/2,-5/6) {$w$};

  \path[every node/.style={font=\sffamily\small}]
   (1) edge node[above]{}(2)
   (2) edge node[above]{}(3)
   (2) edge node[above]{}(5)
   (4) edge node[above]{}(5)
   (5) edge node[above]{}(6);
\end{tikzpicture}
\end{wrapfigure}

Assume $M$ is a complete Riemannian manifold.
A point array $(a_1,\dots,a_n)$ in $M$ together with a choice of a graph with $n$ vertexes labeled by  $(a_1,\dots,a_n)$ and a choice of geodesic $[a_i\,a_j]$ for every adjacent pair $(a_i,a_j)$ is called \emph{geodesic graph}.

For geodesic trees we will use the same notation as for labeled combinatoric tree in square brackets;
for example $[p/xy(q/vw)]$ will denote the geodesic tree with combinatorics as on the diagram.

Fix a geodesic tree $T=[p_1/x_1\dots x_k(p_2/x_{k+1}\dots x_n)]$;
that is, $T$ has two poles $p_1$, $p_2$ and each of the remaining vertexes are adjacent either to $p_1$ or $p_2$ --- the vertexes $x_1,\dots, x_k$ are connected to $p_1$ and $x_{k+1},\dots, x_n$ to $p_2$.

A geodesic tree  $\~T=[\~p_1/\~x_1\dots \~x_k(\~p_2/\~x_{k+1}\dots \~x_n)]$ in the Hilbert space $\HH$ will be called \emph{pivotal tree} for $T$
if 
\begin{enumerate}[(i)]
\item $|\~p_1-\~p_2|_\HH= |p_1-p_2|_M$,
\item $|\~p_i-\~x_j|_\HH= |p_i-x_j|_M$ for any edge $[p_i\,x_j]$ in $T$ and
\item $\measuredangle[\~p_j\,^{\~x_k}_{\~p_i}]_{\HH}=\measuredangle[\~p_j\,^{\~x_k}_{\~p_i}]_M$
for any hinge  $[p_j\,^{x_k}_{\~p_i}]$ in $T$.
\end{enumerate}

\begin{thm}{Rigidity lemma}\label{lem:rigidity}
Let $M$ be a complete Riemannian manifold with nonnegative sectional curvature and $T=[p_1/x_1\dots x_k(p_2/x_{k+1}\dots x_n)]$ be a geodesic tree in $M$.
Suppose  $\~T\zz=[\~p_1/\~x_1\dots \~x_k(\~p_2/\~x_{k+1}\dots \~x_n)]$ is a pivotal tree for  $T$.
Assume that
\[|\~x_i-\~x_j|_\HH\le |x_i-x_j|_M 
\eqlbl{eq:pivotal-comparison}\]
for any pair $(i,j)$ and the convex hull $\~K$ of $\{\~x_1,\dots\~x_n\}$ intersects the line $(\~p_1,\~p_2)$.
Then the equality holds in \ref{eq:pivotal-comparison} for each pair $(i,j)$.
\end{thm}

\parit{Proof.}
Let $\~z$ be a point in the intersection of the line $(\~p_1,\~p_2)$ and $\~K$.
Assume that $\~z$ lies on the half-line from $\~p_1$ to $\~p_2$;
otherwise swap the labels of $\~p_1$ and $\~p_2$.

Denote by $\zeta$ the direction of the geodesic $[p_1\,p_2]$ at $p_1$. 
Set 
\[z=\exp_{p_1}(|\~z-\~p_1|\cdot \zeta).\] 
By comparison, we have
\begin{align*}
|x_i-z|_M &\le |\~x_i-\~z|_{\RR^2}
\end{align*}
for any $i$.

It remains to apply Kirszbraun rigidity theorem (\ref{thm:kirszbraun-rigid}).
\qeds

As above, we assume that $M$ is a complete Riemannian manifold with nonnegative sectional curvature and $[\~p_1/\~x_1\dots \~x_k(\~p_2/\~x_{k+1}\dots \~x_n)]$ is a pivotal tree in $\HH$ for the geodesic tree $[p_1/x_1\dots x_k(p_2/x_{k+1}\dots x_n)]$ in $M$.

Note that by angle comparison, for any $i$ and $j$ we have
\[|\~x_i-\~p_j|_{\HH}\ge |x_i-p_j|_M.\]

It follows that the configuration $\~p_1$, $\~p_2$, $\~x_1,\dots,\~x_n\in \HH$ satisfies the tree comparison (see Section~\ref{sec:intro}) if 
\[|\~x_i-\~x_j|_{\HH}\ge |x_i-x_j|_M
\eqlbl{eq:tree-comparison}\]
for all pairs $(i,j)$.

Let $\Sigma$ be the set of unit vectors normal to the line $(\~p_1, \~p_2)$.
Denote by $\~\xi_i\in \Sigma$ the direction in the half-plane containing $\~x_i$ with the boundary line $(\~p_1, \~p_2)$.

Note that up to a motion of $\HH$, a pivotal configuration is completely described by the angles $\measuredangle(\~\xi_i,\~\xi_j)$.
Moreover, the distance $|\~x_i-\~x_j|_\HH$ is determined by $\measuredangle(\~\xi_i,\~\xi_j)$ and the function $\measuredangle(\~\xi_i,\~\xi_j)\mapsto |\~x_i-\~x_j|_\HH$ is nondecreasing.

Let us denote by $\alpha_{i,j}$ the minimal angle $\measuredangle(\~\xi_i,\~\xi_j)$ in a pivotal configuration such that \ref{eq:tree-comparison} holds.
According to the corollary below $\alpha_{i,j}$ is well defined. 
Note that the inequality \ref{eq:tree-comparison} is equivalent to
\[\measuredangle(\xi_i,\xi_j)\ge \alpha_{i,j}.\]

\begin{thm}{Corollary}\label{cor:|x-x|}
For any geodesic bipolar tree in a complete Riemannian manifold $M$ with nonnegative sectional curvature the following conditions hold:
\begin{enumerate}[(a)]
\item For any pair $i$ and $j$, we have
\[\alpha_{i,j}\le \pi.\]
\item For any triple $i$, $j$ and $k$,  we have
\[\alpha_{i,j}+\alpha_{j,k}+\alpha_{k,i}\le 2\cdot\pi.\]
\end{enumerate}
In other words:
\begin{enumerate}[(a)]
\item\label{cor:|x-x|:a} For any 1(1) geodesic tree (which is a polygonal path) $[p_1/x_1(p_2/x_2)]$ in  $M$ there is a pivotal tree $[\~p_1/\~x_1(\~p_2/\~x_2)]$ such that 
\[|\~x_1-\~x_2|_{\HH}\ge |x_1-x_2|_M.\]

\item\label{cor:|x-x|:b}{\sloppy For any 2(1) geodesic tree $[p_1/x_1x_2(p_2/x_3)]$ in $M$, there is a pivotal tree $[\~p_1/\~x_1\~x_2(\~p_2/\~x_3)]$ such that 
\begin{align*}
|\~x_i-\~x_j|_{\HH}&\ge |x_i-x_j|_M.
\end{align*}
for all $i$ and $j$.

}
\end{enumerate}

\begin{center}
\begin{tikzpicture}[scale=1,
  thick,main node/.style={circle,draw,font=\sffamily\bfseries,minimum size=3mm}]
  \node[main node] (2) at (1,0){$p_1$};
  \node[main node] (3) at (0,0){$x_1$};
  \node[main node] (4) at (3,0) {$x_2$};
  \node[main node] (5) at (2,0) {$p_2$};
  \path[every node/.style={font=\sffamily\small}]
   (2) edge node[above]{}(3)
   (2) edge node[above]{}(5)
   (4) edge node[above]{}(5);
\end{tikzpicture}
\hskip10mm
\begin{tikzpicture}[scale=1,
  thick,main node/.style={circle,draw,font=\sffamily\bfseries,minimum size=3mm}]
  \node[main node] (1) at (0,0) {$x_1$};
  \node[main node] (2) at (1,0){$p_1$};
  \node[main node] (3) at (1,1){$x_2$};
  \node[main node] (4) at (3,0) {$x_3$};
  \node[main node] (5) at (2,0) {$p_2$};

  \path[every node/.style={font=\sffamily\small}]
   (1) edge node[above]{}(2)
   (2) edge node[above]{}(3)
   (2) edge node[above]{}(5)
   (4) edge node[above]{}(5);
\end{tikzpicture}
\end{center}

\end{thm}

\parit{Proof; (\textit{\ref{cor:|x-x|:a}}).}
Consider the pivotal tree $[\~p_1/\~x_1(\~p_2/\~x_2)]$ (which is a polygonal path) with $\measuredangle(\~\xi_1,\~\xi_2)=\pi$.
Note that the points $\~p_1,\~x_1,\~p_2,\~x_2$ are coplanar and the points $\~x_1$ and $\~x_2$ lie on the opposite sides from the line $(\~p_1,\~p_2)$.
It remains to apply the rigidity lemma.

\parit{(\ref{cor:|x-x|:b}).} By (\textit{\ref{cor:|x-x|:a}}), we can assume that \[\alpha_{1,3}+\alpha_{2,3}>\pi.
\eqlbl{sum>pi}\]

Consider the pivotal tree $[\~p_1/\~x_1\~x_2(\~p_2/\~x_3)]$ that lies in a 3-dimensional subspace in such a way that the points $\~x_1$ and $\~x_2$ lie on the opposite sides from the plane containing $\~p_1,\~p_2,\~x_3$, and 
\begin{align*}
\measuredangle(\~\xi_1,\~\xi_3)&=\alpha_{1,3},
&
\measuredangle(\~\xi_2,\~\xi_3)&=\alpha_{2,3}.
\end{align*}
By \ref{sum>pi}, the convex hull $\~K$ of $\{\~x_1,\~x_2,\~x_3\}$ intersects the line $(\~p_1,\~p_2)$.
It remains to apply the rigidity lemma.
\qeds

Note that (\textit{\ref{cor:|x-x|:a}}) and (\textit{\ref{cor:|x-x|:b}}) imply that any nonnegatively curved complete Riemannian manifold satisfies 1(1)-tree and 2(1)-tree comparisons. 
However, 1(1)-tree comparison follows directly from the triangle inequality.

\section{2(2) and 3(1)}\label{6-dipole}

Note that both 2(2)-tree and 3(1)-tree comparisons imply Alexandrov comparison; 
indeed the tripod (that is, 3-tree) is an subtree of both trees 2(2) and 3(1) and the 3-tree comparison is equivalent Alexandrov comparison, see the introduction.
Hence the if part of Theorem~\ref{thm:3(1)+2(2)} follows.

The following proposition is a slightly stronger version of the only-if part --- namely we will show that the required model configuration  in Theorem~\ref{thm:3(1)+2(2)} can be found among pivotal trees.

\begin{center}
\begin{tikzpicture}[scale=1,
  thick,main node/.style={circle,draw,font=\sffamily\bfseries,minimum size=3mm}]

  \node[main node] (1) at (-1/2,0) {$x_1$};
  \node[main node] (2) at (0,1){$p_1$};
  \node[main node] (3) at (-1/2,2){$x_2$};
  \node[main node] (4) at (3/2+.1,0) {$x_4$};
  \node[main node] (5) at (1.1,1) {$p_2$};
  \node[main node] (6) at (3/2+.1,2) {$x_3$};

  \path[every node/.style={font=\sffamily\small}]
   (1) edge node[above]{}(2)
   (2) edge node[above]{}(3)
   (2) edge node[above]{}(5)
   (4) edge node[above]{}(5)
   (5) edge node[above]{}(6);
\end{tikzpicture}
\hskip10mm
\begin{tikzpicture}[scale=1,
  thick,main node/.style={circle,draw,font=\sffamily\bfseries,minimum size=3mm}]

  \node[main node] (1) at (0,0) {$x_1$};
  \node[main node] (2) at (0,1){$p_1$};
  \node[main node] (3) at (-1,1){$x_2$};
  \node[main node] (4) at (2,1) {$x_4$};
  \node[main node] (5) at (1,1) {$p_2$};
  \node[main node] (6) at (0,2) {$x_3$};

  \path[every node/.style={font=\sffamily\small}]
   (1) edge node[above]{}(2)
   (2) edge node[above]{}(3)
   (2) edge node[above]{}(5)
   (2) edge node[above]{}(6)
   (4) edge node[above]{}(5);
\end{tikzpicture}
\end{center}

\begin{thm}{Proposition}\label{2(2)+3(1)}
Let $M$ be a complete Riemannian manifold with nonnegative sectional curvature.
Then for any geodesic 2(2)-tree (or 3(1)-tree) there is a pivotal tree satisfying the corresponding tree comparison.

In particular, $M$ satisfies the 2(2)-tree comparison as well as 3(1)-tree comparison.

\end{thm}

The proofs in the two cases are nearly identical; they differ only by the choice made in the first line.

\parit{Proof.} 
Fix a geodesic tree $[p_1/x_1x_2(p_2/x_3x_4)]$ or $[p_1/x_1x_2x_3(p_2/x_4)]$.
Define the values $\{\alpha_{i,j}\}$ for each pair $i,j$ as in the previous section.

Fix a smooth monotonic function $\phi\:\RR\to\RR$ such that $\phi(x)=0$ if $x\ge 0$ and $\phi(x)>0$ if $x<0$.
Consider a configuration of 4 points $\~\xi_1,\~\xi_2,\~\xi_3,\~\xi_4$ in $\SS^3$ that minimize the \emph{energy}
\[E(\~\xi_1,\~\xi_2,\~\xi_3,\~\xi_4)
=
\sum_{i<j}\phi(\measuredangle(\~\xi_i,\~\xi_j)-\alpha_{i,j}).\]

Consider the geodesic graph $\Gamma$ with 4 vertexes $\~\xi_1,\~\xi_2,\~\xi_3,\~\xi_4$ in $\SS^3$, where 
$\~\xi_i$ is adjacent to $\~\xi_j$ if $\measuredangle(\~\xi_i,\~\xi_j)<\alpha_{i,j}$.
If the comparison does not hold, then $\Gamma$ is not a null graph; that is, $\Gamma$ has some edges.

Note that by Corollary \ref{cor:|x-x|}, degree of any vertex is at least 2.
Indeed, assume $\~\xi_1$ has degree 1 and it is adjacent to $\~\xi_2$, then $\measuredangle(\~\xi_1,\~\xi_2)=\pi$ and therefore $\alpha_{1,2}>\pi$.
The latter contradicts \ref{cor:|x-x|}\textit{\ref{cor:|x-x|:a}}.
If $\~\xi_1$ is a vertex of degree 0, then from the previous case all other vertexes have degree 2.
Since $E$ is minimal the edges $[\~\xi_2\~\xi_3]$, $[\~\xi_3\~\xi_4]$ and $[\~\xi_4\~\xi_1]$ form an equator. 
Therefore $\alpha_{2,3}+\alpha_{2,4}+\alpha_{3,4}>2\cdot \pi$ which contradicts \ref{cor:|x-x|}\textit{\ref{cor:|x-x|:b}}. 

Therefore the graph $\Gamma$ is isomorphic to one the following three graphs.

\begin{center}
\begin{tikzpicture}[scale=1,
  thick,main node/.style={circle,draw,font=\sffamily\bfseries,minimum size=3mm}]

  \node[main node] (1) at (0,0) {};
  \node[main node] (2) at (0,1){};
  \node[main node] (3) at (1,1){};
  \node[main node] (4) at (1,0) {};

  \path[every node/.style={font=\sffamily\small}]
   (1) edge node[above]{}(2)
   (2) edge node[above]{}(3)
   (2) edge node[above]{}(4)
   (3) edge node[above]{}(1)
   (3) edge node[above]{}(4)
   (1) edge node[above]{}(4);
\end{tikzpicture}
\hskip20mm
\begin{tikzpicture}[scale=1,
  thick,main node/.style={circle,draw,font=\sffamily\bfseries,minimum size=3mm}]

  \node[main node] (1) at (0,0) {};
  \node[main node] (2) at (0,1){};
  \node[main node] (3) at (1,1){};
  \node[main node] (4) at (1,0) {};

  \path[every node/.style={font=\sffamily\small}]
   (1) edge node[above]{}(2)
   (2) edge node[above]{}(3)
   (3) edge node[above]{}(1)
   (3) edge node[above]{}(4)
   (1) edge node[above]{}(4);
\end{tikzpicture}
\hskip20mm
\begin{tikzpicture}[scale=1,
  thick,main node/.style={circle,draw,font=\sffamily\bfseries,minimum size=3mm}]

  \node[main node] (1) at (0,0) {};
  \node[main node] (2) at (0,1){};
  \node[main node] (3) at (1,1){};
  \node[main node] (4) at (1,0) {};

  \path[every node/.style={font=\sffamily\small}]
   (1) edge node[above]{}(2)
   (2) edge node[above]{}(3)
   (3) edge node[above]{}(4)
   (1) edge node[above]{}(4);
\end{tikzpicture}
\end{center}

Note that any vertex of $\Gamma$ can not lie in an open hemisphere with all its adjacent vertexes. 
Indeed, if it would be the case, then we could move this vertex increasing the distances to all its adjacent vertexes.
The latter is possible since $\Gamma$ lies in $\SS^3$ (it might be impossible in $\SS^2$).
At the beginning of this move the energy decreases.
Indeed, the portion of the sum in the energy function that corresponds to the edges of $\Gamma$ give linear or quadratic decay, 
the growth of the the sum for the remaining pairs is smaller than quadratic since $\phi(x)=o(x^n)$ for any $n$.
Therefore we arrived to a contradiction.

The 6-edge case (that is, the complete graph with 4 vertexes) can not appear by the rigidity lemma (see \ref{lem:rigidity}).

To do the remaining two cases, note that since the energy is minimal, the angle between the edges at every vertex of degree 2 of $\Gamma$ has to be $\pi$. 
That is, the concatenation of two edges at such vertex is a geodesic.

Consider the 5-edge graph on the diagram.
By the observation above the both triangles in the graph run along one equator.
The latter contradicts Corollary \ref{cor:|x-x|}\textit{\ref{cor:|x-x|:b}}.

For the 4-edge graph (that is, for 4-cycle)
by the same observation the 4 vertexes lie on an equator.%
\footnote{Since any vertex of $\Gamma$ can not lie in an open hemisphere with all its adjacent vertexes, the configuration is centrally symmetric, but we do not need it in the proof.}
Moving one diagonal pair to the north pole 
and the other diagonal pair to the south pole will decrease the energy which leads to a contradiction.
Here we use again that edges condrbute at least quadratic decay of energy, while the nonedges can contribute a slower than quadratic growth.
\qeds

%% file: convexity.tex
\section{Pull-back convexity}\label{convexity}

The part (\textit{\ref{thm:convexity:convexity}}) of Theorem~\ref{thm:convexity} and its converse follows from the three propositions in this section.

\begin{thm}{Proposition}\label{prop:CTIL}
If a complete Riemannian manifold satisfies 4(1)-tree comparison, then it is CTIL.
\end{thm}

\parit{Proof.}
Assume that a Riemannian manifold $M$ satisfies 4(1)-tree comparison.
Assume there is $p\in M$ and $u,v\zz\in \TIL_p$ such that $w=\tfrac12\cdot(u+v)\notin \TIL_p$.
It is sufficient to show that $\gamma(t)=\exp_p(w\cdot t)$ is a length-minimizing on $[0,1]$.

Assume the contrary, that is, $\tau<1$ is the maximal value such that the geodesic $\gamma(t)=\exp_p(w\cdot t)$ is a length-minimizing on $[0,\tau]$.
Set $w'=\tau\cdot w$.
Note that $w'\in\partial \TIL_p$.

Set $q=\exp_p w'$.
By general position argument, we can assume that there are at least two minimizing geodesics connecting $p$ to $q$.\footnote{That is, the set of points that are ends of least two minimizing geodesics from $p$ is dense in the cut locus of $p$.
It was proved in \cite[4.8]{karcher} (the statement is slightly weaker, but the proof proves the needed statement) and latter in \cite{bishop} and \cite{wolter}, see also \cite{petunin-mof} for another proof.}
That is, there is $w''\in \partial \TIL_p$ such that 
$w''\ne w'$
and $\exp_pw'=\exp_pw''$.

\begin{center}
\begin{lpic}[t(-0 mm),b(-0 mm),r(0 mm),l(0 mm)]{pics/7-config(1)}
\lbl[r]{1.5,32;$x$}
\lbl[l]{59.5,31;$y$}
\lbl[t]{18.5,2;$y'$}
\lbl[t]{32.5,2;$x'$}
\lbl[t]{26,4.5;$p$}
\lbl[r]{23,14;$z$}
\lbl[b]{28,31;$q$}
\lbl[l]{30,23;$q'$}
\end{lpic}
\end{center}

Fix small positive real numbers $\delta,\eps$ and $\zeta$.
Consider the following points
\begin{align*}
q'=q'(\eps)&=\exp_p(1-\eps)\cdot w',
&
z=z(\zeta)&=\exp_p(\zeta\cdot w''),
\\
x&=\exp_p u,
&
x'=x'(\delta)&=\exp_p (-\delta\cdot u),
\\
y&=\exp_p v,
&
y'=y'(\delta)&=\exp_p (-\delta\cdot v).
\end{align*}

\begin{wrapfigure}{r}{35 mm}
\begin{tikzpicture}[scale=1,
  thick,main node/.style={circle,draw,font=\sffamily\bfseries,minimum size=8mm}]

  \node[main node] (0) at (3/2,11/6){$x$};
   \node[main node] (1) at (1/2,1/6){$x'$};
  \node[main node] (2) at (3/2,.1){$y$};
  \node[main node] (3) at (1,1){$p$};
  \node[main node] (4) at (1/2,11/6){$y'$};
  \node[main node] (5) at (3,1){$z$};
  \node[main node] (6) at (2,1){$q'$};

  \path[every node/.style={font=\sffamily\small}]
     (0) edge node[above]{}(3)
   (1) edge node[above]{}(3)
   (2) edge node[above]{}(3)
   (3) edge node[above]{}(6)
   (4) edge node[above]{}(3)
   (5) edge node[above]{}(6);
\end{tikzpicture}
\end{wrapfigure}

We will  show that for some choice of $\delta,\eps$ and $\zeta$ the tree comparison for $p/xx'yy'(q'/z)$ does not hold.

Assume the contrary; that is, given any positive numbers $\delta,\eps$ and $\zeta$, there is a point array $\~p$, $\~x$, $\~x'(\delta)$, $\~y$, $\~y'(\delta)$, $\~q'(\eps)$, $\~z(\zeta)\in\HH$ as in the definition of $T$-tree comparison.

If $\delta$ is small, we can assume that $p$ lies on a necessary unique minimizing geodesic $[x\,x']_M$.
Hence 
\[|x-x'|_M=|x-p|_M+|p-x'|_M.\]
By comparison
\begin{align*}
|\~x-\~x'|_\HH&\ge|x-x'|_M,
\\
|\~x-\~p|_\HH&=|x-p|_M,
\\
|\~x'-\~p|_\HH&=|x'-p|_M.
\end{align*}
By triangle inequality,
\[|\~x-\~x'|_\HH=|\~x-\~p|_\HH+|\~x'-\~p|_\HH;\]
that is, $\~p\in [\~x\,\~x']_\HH$.
The same way we see that $\~p\in [\~y\,\~y']_\HH$.

Fix $\eps$ and $\zeta$.
Note that as $\delta\to 0$, we have that  
\begin{align*}
\~x'&\to \~p,
&
\~y'&\to \~p;
\intertext{by first variation formula and tree comparison, we also get that}
\measuredangle[\~p\,^{\~x'}_{\~y}]&\to \measuredangle[p\,^{x'}_{y}],
&
\measuredangle[\~p\,^{\~y'}_{\~x}]&\to \measuredangle[p\,^{y'}_{x}].
\\
\measuredangle[\~p\,^{\~x'}_{\~q'}]&\to \measuredangle[p\,^{x'}_{q'}],
&
\measuredangle[\~p\,^{\~y'}_{\~q'}]&\to \measuredangle[p\,^{y'}_{q'}].
\end{align*}
Indeed, the limits of
$\measuredangle[\~p\,^{\~x'}_{\~y}]$, 
$\measuredangle[\~p\,^{\~x'}_{\~y'}]$,
$\measuredangle[\~p\,^{\~x}_{\~y'}]$
can not be smaller than 
$\measuredangle[p\,^{x'}_{y}]$, 
$\measuredangle[p\,^{x'}_{y'}]$,
$\measuredangle[p\,^{x}_{y'}]$
respectfully (the last three angles do not depend on $\delta$).
On the other hand, the limits of $\measuredangle[\~p\,^{\~x'}_{\~y}]+\measuredangle[\~p\,^{\~x'}_{\~y'}]$ and $\measuredangle[\~p\,^{\~x}_{\~y'}]+\measuredangle[\~p\,^{\~x'}_{\~y'}]$ can not be bigger than $\pi$.
Since $\measuredangle[p\,^{x'}_{y}]+\measuredangle[p\,^{x'}_{y'}]=\pi$ and $\measuredangle[p\,^{x}_{y'}]+\measuredangle[p\,^{x'}_{y'}]=\pi$ the first two identities follow.
The last two identities follow the same ways since \[\measuredangle[\~p\,^{\~x'}_{\~q}]+\measuredangle[\~p\,^{\~y'}_{\~q}]+\measuredangle[\~p\,^{\~x'}_{\~y'}]\le 2\cdot \pi\]
and 
\[\measuredangle[p\,^{x'}_{q}]+\measuredangle[p\,^{y'}_{q}]+\measuredangle[p\,^{x'}_{y'}]= 2\cdot \pi.\]

Therefore
\begin{align*}
\measuredangle[\~p\,^{\~x}_{\~y}]&\to \measuredangle[p\,^x_y],
&
\measuredangle[\~p\,^{\~x}_{\~q'}]&\to \measuredangle[p\,^x_{q'}],
&
\measuredangle[\~p\,^{\~y}_{\~q'}]&\to \measuredangle[p\,^y_{q'}].
\end{align*}

Therefore, passing to a partial limit as $\delta\to0$, we get a configuration of 5 points 
$\~p, \~x,\~y,\~q'=\~q'(\eps),\~z=\~z(\zeta)$ such that  
\begin{align*}
\measuredangle[\~p\,^{\~x}_{\~y}]&= \measuredangle[p\,^{x}_{y}],
&
\measuredangle[\~p\,^{\~y}_{\~q'}]&= \measuredangle[p\,^{y}_{q'}],
&
\measuredangle[\~p\,^{\~x}_{\~q'}]&= \measuredangle[p\,^{x}_{q'}].
\end{align*}
In other words, the map sending the points $0,u,v,w'\in\T_p$ to $\~p,\~x,\~y,\~q'\in\HH$ correspondingly is distance preserving.

Note that $q'\to q$ as $\eps\to0$. 
Therefore, in the limit,
we get a configuration $\~p$, $\~x$, $\~y$, $\~q'$, $\~z=\~z(\zeta)$ such that in addition we have
\begin{align*}
|\~q'-\~z|&=|q-z|,
&
|\~p-\~z|&\ge |p-z|,
\\
|\~x-\~z|&\ge |x-z|,
&
|\~y-\~z|&\ge |y-z|.
\end{align*}

Since $w''\ne w'$, for small values $\zeta$ the last three inequalities 
imply 
\[|\~q'-\~z|>|q-z|,\]
a contradiction.
\qeds

\begin{thm}{Proposition}\label{prop:convex}
If  a complete CTIL Riemannian manifold $M$ satisfies 4(1)-tree comparison,
then for any $p,q\in M$, we have $f''\le 1$, where $f$ is the function $f\: \TIL_p\to \RR$ defined by
\[f(v)=\tfrac12\cdot\dist_q^2\circ\exp_p(v).\]

\end{thm}

\parit{Proof.}
Note that 4(1)-tree comparison implies 3-tree comparison.
Hence $M$ has nonnegative sectional curvature.

Fix $u,v\in \TIL_p$ and $w\in [u\,v]$.
It is sufficient to show that there is a function $g\:\T_p\to \RR$ such that
\[g''=1,\quad
g(w)=f(w),\quad
g(u)\ge f(u)\quad
\text{and}\quad
g(v)\ge f(v).\]

Fix small $\eps>0$ and set
\begin{align*}
x&=\exp_p u,
&y&=\exp_p v, 
&z&=\exp_pw,
\\
x'&=\exp_p(-\eps\cdot  u),
&y'&=\exp_p(-\eps\cdot  v).
\end{align*}
Let us apply the $p/xyx'y'(z/q)$ comparison and pass to the limit as $\eps\to 0$
as we did in the proof of Proposition~\ref{prop:CTIL}.
We obtain a configuration of points $\~p, \~x, \~y, \~z, \~q\in\HH$, satisfying corresponding comparisons and
in addition
\begin{align*}
\measuredangle[\~p\,^{\~x}_{\~y}]&= \measuredangle[p\,^{x}_{y}],
&\measuredangle[\~p\,^{\~x}_{\~z}]&= \measuredangle[p\,^{x}_{z}],
&\measuredangle[\~p\,^{\~z}_{\~y}]&= \measuredangle[p\,^{z}_{y}].
\end{align*}
In particular,
from above and Toponogov comparison, we have
\begin{align*}
|\~x-\~y|_\HH&=|u-v|_{\T_p},
&|\~z-\~y|_\HH&=|w-v|_{\T_p},
&|\~x-\~z|_\HH&=|u-w|_{\T_p},
\\
|\~q-\~z|_\HH&=|q-z|_M,
&|\~q-\~x|_\HH&\ge|q-x|_M,
&|\~q-\~y|_\HH&\ge|q-y|_M.
\end{align*}
In particular, there is a distance-preserving map $\T_p\to \HH$ 
such that $u\mapsto \~x$, $v\mapsto \~y$, $w\mapsto \~z$ and $0\mapsto \~p$.
Further, we identify $\T_p$ and a subset of $\HH$ using this map.

Consider the function $g(s):=\tfrac12\cdot|s-\~q|_{\T_p}^2$.
Note that $g''=1$ and
\begin{align*}
g(w)&=\tfrac12\cdot|\~q-\~z|_{\T_p}^2
=\tfrac12\cdot|q-z|_M^2
=f(w),
\\
g(u)&=\tfrac12\cdot|\~q-\~x|_{\T_p}^2
\ge \tfrac12\cdot|q-x|_M^2
=f(u),
\\
g(v)&=\tfrac12\cdot|\~q-\~y|_{\T_p}^2
\ge\tfrac12\cdot|q-y|_M^2
= f(u).
\end{align*}
Hence the statement.
\qeds

\begin{thm}{Proposition}\label{prop:m(n)}
Assume $M$ is a complete CTIL Riemannian manifold such that for any $p,q\in M$, we have $f''\le 1$, where $f$ is the function $f\: \TIL_p\to \RR$ defined by
\[f(v)=\tfrac12\cdot\dist_q^2\circ\exp_p(v),\]
Then $M$ satisfies all bipolar comparisons.
\end{thm}

\parit{Proof.}
Fix points $p$ and $q$ in $M$;
set $\~q=\log_pq\in\T_p$ and $\~f(v)=\tfrac12\cdot |v-\~q|_{\T_p}^2$.
Note that $d_{\~q}f=d_{\~q}\~f$ and $\~f''\equiv 1$, therefore
\[f\le \~f.\eqlbl{eq:f=<f}\]
Further note that the inequality \ref{eq:f=<f} is equivalent to the Toponogov comparison for all hinges $[p\,{}^x_q]$ in $M$.
It follows that $M$ has nonnegative sectional curvature. 

\medskip

Fix a bipolar geodesic tree $[p/x_1\dots x_n(q/y_1\dots y_m)]$ in $M$.
Set 
\[\~p=0=\log_pp,\quad \~q=\log_pq,\quad\text{and}\quad \~x_i=\log_px_i\]
for each $i$. 

Let $\psi_1\:\T_q\to \T_p$ be the dual of the differential $d_{\~q}\exp_p\:\T_p=\T_{\~q}\T_p\to \T_q$; note that for any smooth function $h$
\[\psi_1\:\nabla_{q}h\mapsto \nabla_{\~q}(h\circ\exp_p).\]
Since sectional curvature of $M$ is nonnegative, the restriction $\exp_p|_{\TIL_p}$ is short and therefore so is $\psi_1$.

In particular there is a linear map $\psi_2\:\T_q\to\T_p$ such that, the map $\iota\:\T_q\zz\to \T_p\oplus\T_p$ defined by
\[\iota\:v\mapsto (\psi_1(v), \psi_2(v))\]
is distance preserving.

Further set 
\begin{align*}
h_i&=\tfrac12\cdot\dist_{y_i}^2,
&
g_i&=h_i\circ\exp_p|_{\TIL_p},
\\
\~y_i&=\~q-\psi_1(\nabla_q h_i),
&
\~z_i&=-\psi_2(\nabla_q h_i).
\end{align*}

By construction
\[|(\~y_i,\~z_i)-(\~q,0)|_{\T_p\oplus\T_p}=|y_i-q|_M.\]

At the point $(\~q,0)\in \T_p\oplus\T_p$ the restriction functions $\~g_i=\tfrac12\cdot\dist^2_{\~y_i}|_{\T_p\oplus 0}$ and the function $g_i$ have the same value and gradient.
Since $g_i''\le 1$ and $\~g_i''=1$, we get $\~g_i\ge g_i$. 
The latter implies
\begin{align*}
|(\~y_i,\~z_i)-(\~p,0)|_{\T_p\oplus\T_p}&\ge|y_i-p|_M
\\
|(\~y_i,\~z_i)-(\~x_j,0)|_{\T_p\oplus\T_p}&\ge|y_i-x_j|_M.
\end{align*}
for any $i$ and $j$.
That is, the configuration 
\[(\~p,0), (\~x_1,0),\dots,(\~x_n,0), (\~q,0), (\~y_1,\~z_1),\dots,(\~y_m,\~z_m)\]
satisfies the comparison.
It remains to apply to it an isometric embedding $\T_p\oplus\T_p\hookrightarrow \HH$.
\qeds

%% file: mtw.tex
\section{MTW}\label{MTW+}

Proposition~\ref{MTW-plus-convexity} provides equivalence of (\textit{\ref{thm:convexity:convexity}}) and (\textit{\ref{thm:convexity:MTW}}) in the main theorem (\ref{thm:convexity}); which is the final step its proof.
The equivalence is proved by calculations along the same lines as in \cite[Chapter 12]{villani}.

Let us introduce notations and use them to reformulate the property (\textit{\ref{thm:convexity:MTW}}).

\parbf{Tangent vectors.}
Let $M$ be a Riemannian manifold, $p\in M$.
Denote by $\IL_p$ the \emph{inner locus} of $p$; it can be defined as the $\exp_p$-image of $\TIL_p$ or, equivalently, as the complement $M\backslash \CL_p$, where $\CL_p$ denotes the cut locus of $p$.
Note that $q\in\IL_p$ if and only if $p\in \IL_p$.

Assume $q\in\IL_p$; that is, $q=\exp_pW$ for some $W\zz\in \TIL_p M$.
Given a vector $Y\in T_q$, consider the unique vector $Y_p\in\T_p$ such that 
\[Y=(d_W\exp_p)Y_p.\]

Note that $p=\exp_q(-W_q)$ if $p$, $q$ and $W$ are as above.

Given $x\in\IL_p$ such that $x=\exp_pX$ for some  $X\in \TIL_p$,  set
\[\tilde Y_p(x)=(d_X\exp_p) Y_p;\]
this way we defined a vector field $\tilde Y_p$ in $\IL_p$.

Note that in the vector field $\tilde Y_p$ is constant in the normal coordinates at $p$;
in particular 
\[\nabla_X\tilde Y_p=0\eqlbl{eq:zero}\] 
for any $X\in\T_p$.
Further, note that 
\[Y\tilde Y_pf=Y\tilde Y_qf+(\nabla_Y\tilde Y_p)f\eqlbl{eq:nabla}\]
for $Y\in \T_q$ and any smooth function $f$.
Indeed applying \ref{eq:zero}, we get that
\begin{align*}
(Y\tilde Y_p-Y\tilde Y_q)f
&=[(\tilde Y_q\tilde Y_p-\tilde Y_p\tilde Y_q)f](q)=
\\
&=(\nabla_Y \tilde Y_p-\nabla_Y\tilde Y_q)f=\nabla_Y \tilde Y_pf.
\end{align*}

\parbf{Column notation.} Given two points $p$ and $q$ in a Riemannian manifold $M$,
let us define the cost function $(p,q)\mapsto \cost{p}{q}$ as
\[\dcost{}{}{p}{q}=\tfrac12\cdot|p-q|^2_M\]
We will need to differentiate the cost function by both argument.
In order to avoid possible confusion, we will write the vector next to the differentiated argument. 
For example
\[\dcost{X}{Y}{p}{q}\]
is the second mixed derivative of the cost function at the pair $(p,q)$, once by the first argument ($p$) along the vector $X \in \T_p$ and once by the second argument ($q$) along the vector field $Y\in \T_q$.
We may also write a vector field instead of the vectors.

Using the introduced notations,
we can reformulate the property (\textit{\ref{thm:convexity:MTW}}) in Theorem~\ref{thm:convexity}
as

\begin{itemize}
 \item[\textit{(ii)}$'$] \emph{If $X\in \T_p$, $Y\in\T_q$ and $q\in \IL_p$, then
 \[\dcost{X\tilde X_p}{Y\tilde Y_p}{p}{q}\le 0.\]}
\end{itemize}

The left hand side of the last inequality, multiplied by $(-\tfrac32)$ is called \emph{MTW-curvature} or \emph{cost-curvature}; it is denoted by $\mathfrak{S}(X,Y)$, see \cite[equation 12.21]{villani};
if $p=q$, then $\mathfrak{S}(X,Y)$ coincides with the curvature $\langle\Rm(X,Y)Y,X\rangle$, see \cite[12.30]{villani}.
In particular, if the condition (\textit{\ref{thm:convexity:MTW}}) holds, then the manifold has nonnegative sectional curvature.

Assume $q=\exp_pW$ for some $W\in\TIL_p$. 
Then 
\[\dcost{X}{{}}{p}{q}=-\langle X,W\rangle;\eqlbl{derX}\]
\[\dcost{X}{Y}{p}{q}=-\langle X,Y_p\rangle =-\langle X_q,Y\rangle\eqlbl{derXY}\]
and
\[\dcost{X}{Y\tilde Y_p}{p}{q}=0.\eqlbl{derXYY}\]

Indeed, \ref{derX} is equivalent to the first variation formula.
Taking the derivative of \ref{derX} in the normal coordinates at $p$ we get \ref{derXY}.
\[\dcost{X}{Y}{p}{q}=-\langle X,Y_p\rangle.\]
Since $q\in \IL_p$ if and only if $p\in\IL_q$ we can swap $p$ and $q$ and get the second identity in \ref{derXY}. Finally, the value $\langle X,Y_p\rangle$ does not depend on $q$; therefore the derivative along the second argument \ref{derXY} has to vanish;
hence \ref{derXYY} follows.

Let us use the identities to show that 
\[\dcost{X\tilde X_p}{Y\tilde Y_p}{p}{q}=\dcost{X\tilde X_p}{Y\tilde Y_p}{p}{q}
\quad
\text{or, equivalently}
\quad
\mathfrak{S}(X,Y)=\mathfrak{S}(Y,X).\eqlbl{derXXYY}\]
This identity will not be used in the sequel, but it might help the reader to adapt to the column notation.

Applying \ref{eq:nabla}, we get that
\begin{align*}
\dcost{X\tilde X_p}{Y\tilde Y_p}{p}{q}&=\dcost{X\tilde X_p}{Y\tilde Y_q}{p}{q}+\dcost{X\tilde X_p}{\nabla_Y\tilde Y_q}{p}{q}=
\\
&=\dcost{X\tilde X_p}{Y\tilde Y_q}{p}{q}+\dcost{X\tilde X_q}{\nabla_Y\tilde Y_q}{p}{q}-\dcost{\nabla_X\tilde X_p}{\nabla_Y\tilde Y_q}{p}{q}.
\end{align*}
By \ref{derXYY}, 
\[\dcost{X\tilde X_q}{\nabla_Y\tilde Y_q}{p}{q}=0.\]
Therefore
\[\dcost{X\tilde X_p}{Y\tilde Y_p}{p}{q}=\dcost{X\tilde X_p}{Y\tilde Y_q}{p}{q}-\dcost{\nabla_X\tilde X_p}{\nabla_Y\tilde Y_q}{p}{q}.\]
The right hand side is symmetric in $p$ and $q$;
hence \ref{derXXYY} follows.

\begin{thm}{Proposition}\label{MTW-plus-convexity}
Let $M$ be a CTIL Riemannian manifold.
Then the following conditions are equivalent:
\begin{enumerate}[(a)]
 \item\label{MTW-plus-convexity:MTW} For any $p\in M$, $q\in \IL_p$, $X\in \T_p$ and $Y\in\T_q$ we have
 \[\dcost{X\tilde X_p}{Y\tilde Y_p}{p}{q}\le 0.\]
 \item\label{MTW-plus-convexity:h} For any $p_0,p_1\in M$, the function $h\:\TIL_{p_0}\to\RR$ defined by
\[h(X)=\cost{{p_1}}{\exp_{p_0}X}-\cost{{p_0}}{\exp_{p_0}X}\]
is concave;
 \item\label{MTW-plus-convexity:f}  For any $p_0,p_1\in M$, the function $f\:\TIL_{p_0}\to\RR$ defined by
\[f(X)=\cost{p_1}{\exp_{p_0}X}\]
is 1-concave.
\end{enumerate}
\end{thm}

Note that equivalence 
(\textit{\ref{MTW-plus-convexity:MTW}})$\iff$(\textit{\ref{MTW-plus-convexity:f}}) imply the
equivalence of (\textit{\ref{thm:convexity:convexity}})$\iff$(\textit{\ref{thm:convexity:MTW}}) in \ref{thm:convexity}. 
Therefore the proposition finishes the proof of the main theorem.

\parit{Proof.}  
Note that $\cost{p_0}{\exp_{p_0}X}=\tfrac12\cdot |X|^2$ for any $X\in\TIL_{p_0}$,
in particular the function $s(X)\zz=\cost{p_0}{\exp_{p_0}X}$ is $1$-affine (that is, $1$-concave and $1$-convex at the same time).

Evidently $f=h+s$, therefore (\textit{\ref{MTW-plus-convexity:h}})$\iff$(\textit{\ref{MTW-plus-convexity:f}}).

\parit{(\ref{MTW-plus-convexity:MTW}) $\Rightarrow$ (\ref{MTW-plus-convexity:h}).}
Note that the function $h$ is semiconcave.
Therefore $h$ is concave if
\[\tfrac{d^2}{dt^2}h(U+t\cdot V)\le 0\]
at $t=0$ for almost all vectors $U\in\TIL_p$ and $V\in \T_p$.

Using the column notation, we can rewrite the inequality in the following equivalent form:
\[
\begin{matrix}{{}}\\{Y\tilde Y_{p_0}}
\end{matrix}
\left(\dcost{}{}{p_1}{q}-\dcost{}{}{p_0}{q}\right)\le 0\eqlbl{eq:MTW-plus}\]
for any $p_0,p_1, q$ and $Y\in \T_q$;
from above it is sufficient to prove \ref{eq:MTW-plus} for almost all $q$; in particular, we can assume that $p_0,p_1\in \IL_q$.

Let $W, X\in \TIL_q$ be such that $p_0=\exp_qW$, $p_1=\exp_q(W+X)$.
Since $M$ is CTIL, $W+t\cdot X\in\TIL_q$ for any $t\in[0,1]$;
set $p_t=\exp_q(W+t\cdot X)$.

Let us use the identity
$f(1)-f(0)-f'(0)=\int_0^1f''(t)\cdot(1-t)\cdot dt$,
for the function 
\[f(t)=\dcost{}{Y\tilde Y_{q}}{p_t}{q};\]
Note that
\[f'(0)=\dcost{X}{Y\tilde Y_{q}}{p_0}{q}
\quad
\text{and}
\quad 
f''(t)
=\dcost{\tilde X_{q}\tilde X_{q}}{Y\tilde Y_{q}}{p_t}{q},\]
therefore
\[\dcost{}{Y\tilde Y_{q}}{p_1}{q}-\dcost{}{Y\tilde Y_{q}}{p_0}{q}-\dcost{X}{Y\tilde Y_{q}}{p_0}{q}=\int_0^1 \dcost{\tilde X_{q}\tilde X_{q}}{Y\tilde Y_{q}}{p_t}{q}\cdot dt.\]

By (\textit{\ref{MTW-plus-convexity:MTW}}), the term under the integral is nonpositive; therefore
\[
\dcost{{}}{Y\tilde Y_{q}}{}{}
\left(\dcost{}{}{p_1}{q}-\dcost{}{}{p_0}{q}\right)
\le
\dcost{X}{Y\tilde Y_{q}}{p_0}{q}.\]
By \ref{eq:nabla}, we can rewrite the last inequality the following way:
\[
\dcost{{}}{Y\tilde Y_{p_0}}{}{}
\left(\dcost{}{}{p_1}{q}-\dcost{}{}{p_0}{q}\right)
\le
\dcost{{}}{\nabla_Y \tilde Y_{p_0}}{}{}
\left(\dcost{}{}{p_1}{q}
-
\dcost{}{}{p_0}{q}\right)
+
\dcost{X}{Y\tilde Y_{q}}{p_0}{q}.
\eqlbl{eq:MTW-plus-extra}\]

Applying \ref{derX}, we get that
\begin{align*}
\dcost{{}}{\nabla_Y \tilde Y_{p_0}}{}{}
\left(\dcost{}{}{p_1}{q}
-
\dcost{}{}{p_0}{q}\right)
&=
-\langle\nabla_Y \tilde Y_{p_0},W+X\rangle + \langle\nabla_Y \tilde Y_{p_0},W\rangle
=
\\
&=-\langle\nabla_Y \tilde Y_{p_0},X\rangle.
\end{align*}
Further, applying \ref{derXYY} and \ref{derXY}, we get that
\begin{align*}
\dcost{X}{Y\tilde Y_{q}}{p_0}{q}
&=
\dcost{X}{Y\tilde Y_{p_0}}{p_0}{q}-\dcost{X}{\nabla_Y\tilde Y_{p_0}}{p_0}{q}=
\\
&=0+\langle\nabla_Y \tilde Y_{p_0},X\rangle.
\end{align*}
It follows that the right hand side in \ref{eq:MTW-plus-extra} vanishes;
hence \ref{eq:MTW-plus} follows.

\parit{(\ref{MTW-plus-convexity:h}) $\Rightarrow$ (\ref{MTW-plus-convexity:MTW}).}
Let $p_t$, $q$, $W$, $X$ and $Y$ be as above;
set 
\[h_t(Z)=\cost{{p_t}}{\exp_{p_0}Z}-\cost{{p_0}}{\exp_{p_0}Z}.\]
Note that $h_0\equiv 0$; in particular 
\[Y\tilde Y_{p_0}h_0=0\]
for any $Y\in\T_q$.

By (\textit{\ref{MTW-plus-convexity:h}}), 
\[Y\tilde Y_{p_0}h_t\le 0.\]
It follows that 
\[\frac{d^2}{dt^2}(Y\tilde Y_{p_0}h_t)\le 0\]
at $t=0$.
Finally note that 
\[\dcost{X\tilde X_{p_0}}{Y\tilde Y_{p_0}}{p_0}{q}=\frac{d^2}{dt^2}(Y\tilde Y_{p_0}h_t);\]
hence the part (\textit{\ref{MTW-plus-convexity:h}}) follows.
\qeds

%% file: polypole.tex
\section{All tree comparisons}\label{sec:all-tree}

\parit{Proof of Theorem~\ref{thm:hilbert-quotient}.}
The ``if'' part follows from Exercise~\ref{ex:quotient};
let us prove the ``only if'' part.

Since $X$ is separable, it has a countable everywhere dense set $\{x_1,x_2,\dots\}$.
Consider the complete countable graph $K$ with the vertexes labeled by $x_1,x_2\dots$;
let $T\to K$ be its universal covering.

The graph $T$ is a tree with countable set of vertexes.
Note that $K$ can be presented as a union of a nested sequence of finite trees $T_1\subset T_2\subset \dots$;
moreover, we can assume that each $T_n$ is spanned by vertexes $\{y_1,\dots,y_n\}$ for some enumeration $y_1,y_2,\dots$ of the vertexes of $T$.
Denote by $s(y_i)$ the corresponding point $x_j$ in $X$.

Applying the tree comparison for each $T_n$ we get a finite configuration $\~y_{1,n},\dots,\~y_{n,n}$ in $\HH$.
Set $Y_n=\{\~y_{1,n},\dots,\~y_{n,n}\}$. 
Consider the map $s_n\:Y_n\to X$ defined by $s_n\:\~y_{i,n}\mapsto s(y_i)$.
By comparison, $s_n$ is a short map that preserves the distances between the points adjacent in $T_n$.

Without loss of generality, we may assume that for any $k\le n$, the points $\~y_{1,n},\dots,\~y_{k,n}$ lie in the subspace spanned by first $k-1$ elements of a fixed basis of $\HH$.
Passing to a partial limit as $n\to \infty$, we get a set $Y=\{\~y_1,\~y_2,\dots\}\subset \HH$ and a map $s\:Y\to X$ that is short and preserves the distances between the points adjacent in $T$.
In particular for any $\~y_i$ and $x_j$ there is $\~y_k$ such that $s(\~y_k)=x_j$ and $|\~y_i-\~y_k|_\HH=|s(\~y_i)-x_j|_X$.

By construction, the continuous extension of $s$ to the closure $\bar Y$ of $Y$ is a required submetry.
\qeds

The following proof was suggested by Alexander Lytchak, another proof follows from the construction of Chuu-Lian Terng and Gudlaugur Thorbergsson given in \cite[Section 4]{terng-thorbergsson}.

\parit{Proof of Proposition~\ref{prop:group}.}
Denote by $G^n$ the direct product of $n$ copies of $G$.
Consider the map $\phi_n\:G^n\to G/\!\!/H$ defined by
\[\phi_n\:(\alpha_1,\dots,\alpha_n)\mapsto [\alpha_1\cdots\alpha_n]_H,\]
where $[x]_H$ denotes the $H$-orbit of $x$ in $G$.

Note that $\phi_n$ is a quotient map for the action of $H\times G^{n-1}$ on $G^n$ defined by
\[(\beta_0,\dots,\beta_n)\cdot(\alpha_1,\dots,\alpha_n)=(\gamma_1\cdot \alpha_1\cdot\beta_1^{-1},\beta_1\cdot\alpha_2\cdot\beta_2^{-1},\dots,\beta_{n-1}\cdot\alpha_n\cdot\beta_n^{-1}),\]
where $\beta_i\in G$ and $(\beta_0,\beta_n)\in H<G\times G$.

Denote by $\rho_n$ the product metric on $G^n$ rescaled with factor $\sqrt{n}$.
Note that the quotient $(G^n,\rho_n)/(H\times G^{n-1})$ is isometric to $G/\!\!/H=(G,\rho_1)/\!\!/H$.
Let $\phi_n\:(G^n,\rho_n)\to G/\!\!/H$ be the corresponding quotient map; 
clearly $\phi_n$ is a submetry. 

As $n\to\infty$ the curvature of $(G^n,\rho_n)$ converges to zero and its injectivity radius goes to infinity.
Therefore the ultra-limit of $(G^n,\rho_n)$ with marked identity element is a Hilbert space $\HH$ and the submetries $\phi_n$ ultra-converge to a submetry $\phi\:\HH\to G/\!\!/H$.
It remains to apply Exercise~\ref{ex:quotient}.
\qeds

%% file: end.tex
\section{Remarks}

\parbf{On graph comparison.}
One can define \emph{graph comparison} for any graph by stating that there is a model configuration in $\HH$ such that 
\begin{enumerate}[(i)]
\item\label{ge} the distance between each pair of adjacent points is at most as big 
and
\item the distance between each pair of nonadjacent is at least as big.
\end{enumerate}

Note that in the definition of tree comparison we used equalities instead of inequalities in (\ref{ge}). 
However the new definition agrees with the old one for  trees: 

\begin{thm}{Exercise}
Show that if a graph is a tree, then the graph comparison defined above is equivalent to the tree comparison defined at the beginning of the paper.
\end{thm}

Note that nonnegative and nonpositive curvature can be defined using the comparison for following two graphs on 4 vertexes:

\begin{center}
\begin{tikzpicture}[scale=1,
  thick,main node/.style={circle,draw,font=\sffamily\bfseries,minimum size=3mm}]

  \node[main node] (0) at (0,0) {};
  \node[main node] (1) at (-3/5,1/3){};
  \node[main node] (2) at (3/5,1/3){};
  \node[main node] (3) at (0,-2/3) {};

  \path[every node/.style={font=\sffamily\small}]
   (0) edge node[above]{}(1)
   (0) edge node[above]{}(2)
   (0) edge node[above]{}(3);
\end{tikzpicture}
\hskip30mm
\begin{tikzpicture}[scale=1,
  thick,main node/.style={circle,draw,font=\sffamily\bfseries,minimum size=3mm}]

  \node[main node] (1) at (0,0) {};
  \node[main node] (2) at (0,1){};
  \node[main node] (3) at (1,1){};
  \node[main node] (4) at (1,0) {};

  \path[every node/.style={font=\sffamily\small}]
   (1) edge node[above]{}(2)
   (1) edge node[above]{}(4)
   (2) edge node[above]{}(3)
   (3) edge node[above]{}(4);
\end{tikzpicture}
\end{center}
If a graph $G$ has two induced subgraphs that isomorphic to each of these two graphs, then the corresponding graph comparison implies that the curvature vanish in the sense of Alexandrov.
In particular, any complete length spaces satisfying $G$-graph comparison is isometric to a convex set in a Hilbert space. 

By Reshetnyak majorization theorem,
the nonpositive curvature could be also defined using the comparison for cycle;
for example the 6-cycle --- the first graph on the following diagram.

\begin{center}
\begin{tikzpicture}[scale=1,
  thick,main node/.style={circle,draw,font=\sffamily\bfseries,minimum size=3mm}]

  \node[main node] (4) at (3/2,5/6){};
   \node[main node] (1) at (1/2,-5/6){};
  \node[main node] (2) at (3/2,-5/6){};
  \node[main node] (0) at (0,0){};
  \node[main node] (5) at (1/2,5/6){};
  \node[main node] (3) at (2,0){};

  \path[every node/.style={font=\sffamily\small}]
    (0) edge node[above]{}(1)
   (1) edge node[above]{}(2)
    (2) edge node[above]{}(3)
   (3) edge node[above]{}(4)
   (4) edge node[above]{}(5)
   (5) edge node[above]{}(0);
\end{tikzpicture}
\hskip30mm
\begin{tikzpicture}[scale=1,
  thick,main node/.style={circle,draw,font=\sffamily\bfseries,minimum size=3mm}]

  \node[main node] (4) at (3/2,5/6){};
   \node[main node] (1) at (1/2,-5/6){};
  \node[main node] (2) at (3/2,-5/6){};
  \node[main node] (0) at (0,0){};
  \node[main node] (5) at (1/2,5/6){};
  \node[main node] (3) at (2,0){};

  \path[every node/.style={font=\sffamily\small}]
    (0) edge node[above]{}(1)
   (1) edge node[above]{}(2)
    (2) edge node[above]{}(3)
   (3) edge node[above]{}(4)
   (4) edge node[above]{}(5)
   (5) edge node[above]{}(0)
    (0) edge node[above]{}(2)
   (1) edge node[above]{}(3)
    (2) edge node[above]{}(4)
   (3) edge node[above]{}(5)
   (4) edge node[above]{}(0)
   (5) edge node[above]{}(1);
\end{tikzpicture}
\end{center}

The comparison for the octahedron graph (the second graph on the diagram) implies that the space is nonpositively curved.
The latter follows since in this graph, a 4-cycle appears as an induced subgraph.
On the other hand, this comparison might be stronger and it might be interesting to understand.
(Note that the octahedron graph has no induced tripods.) 

The quotients of Hilbert space provide motivating examples to consider the tree comparison (see Theorem~\ref{thm:hilbert-quotient}).
Unfortunately we do not have such guiding examples in nonpositive curvature --- we can only fumble around for something without light.

\parbf{On colored graph comparison.}
It is also possible to use a graph with $(\mp)$-colored edges and define comparison by model configuration such that the distances between vertexes adjacent by a $(-)$-edge do not get larger and by $(+)$-edge do not get smaller and no condition on the remaining pairs.

\begin{center}
\begin{tikzpicture}[scale=1,
  thick,main node/.style={circle,draw,font=\sffamily\bfseries,minimum size=3mm}]

  \node[main node] (0) at (0.5,0){};
   \node[main node] (1) at (1,1){};
  \node[main node] (2) at (1,-1){};
  \node[main node] (3) at (2,1){};
  \node[main node] (4) at (2,-1){};
  \node[main node] (5) at (3,1){};
\node[main node] (6) at (3,-1){};
 \node[main node] (7) at (4,1){};
\node[main node] (8) at (4,-1){};
\node[main node] (100) at (4.5,0){};

  \path[every node/.style={font=\sffamily\small}]
    (0) edge[dashed] node[above]{}(100)
   (1) edge[dashed] node[above]{}(2)
    (3) edge[dashed] node[above]{}(4)
   (5) edge[dashed] node[above]{}(6)
   (7) edge[dashed] node[above]{}(8)
   (0) edge node[above]{}(1)
   (0) edge node[above]{}(2)
    (1) edge node[above]{}(3)
   (2) edge node[above]{}(4)
   (1) edge node[above]{}(4)
   (2) edge node[above]{}(3)
   (5) edge node[above]{}(3)
   (6) edge node[above]{}(4)
   (5) edge node[above]{}(4)
   (6) edge node[above]{}(3)
   (5) edge node[above]{}(7)
   (6) edge node[above]{}(8)
   (5) edge node[above]{}(8)
   (6) edge node[above]{}(7)
   (7) edge node[above]{}(100)
   (8) edge node[above]{}(100);
\end{tikzpicture}
\end{center}

For example, the $(2{\cdot}n+2)$-comparison (which holds in $\CAT[0]$ length spaces, see \cite{AKP}) can be considered as a comparison for the colored graph above, where $(-)$-edges are marked by solid lines and $(+)$-edges by dashed lines.

\parbf{Finite subsets of Alexandrov spaces.}
The following problem discussed in \cite[7.1]{AKP} was one of the original motivations to study the tree comparison:
\emph{Which finite metric spaces admit isometric embeddings into some Alexandrov spaces with nonnegative curvature.}

This problem is still open (as well as its analog for $\CAT(0)$ spaces).
The $(n-1)$-tree comparison provides a necessary condition for $n$-point metric spaces, see \cite[4.1]{AKP}.
(The so called \emph{matrix inequality} turns out to be weaker; see the discussion below.)
This condition is sufficient for the 4-point metric spaces.
It might be still sufficient for 5-point metric spaces,
but not for 6-point metric spaces.

The corresponding example of 6-point metric space was constructed by Sergei Ivanov, see \cite{AKP}.
Theorem~\ref{2(2)+3(1)}, provides a source for such examples --- any 6-point metric space that satisfy all 5-tree comparisons, but does not satisfy 2(2)-tree comparison provides an example.
This class of examples includes the example of Sergei Ivanov --- in the notations of \cite[7.1]{AKP} it does not satisfies the comparison for the tree $y/az(q/xb)$.

By Theorem~\ref{2(2)+3(1)} and a theorem in \cite{AKP}, 
5-tree and 2(2)-tree comparisons provide a necessary condition for 6-point metric spaces.
We expect that these conditions are sufficient

\begin{center}
\begin{tikzpicture}[scale=1,
  thick,main node/.style={circle,draw,font=\sffamily\bfseries,minimum size=3mm}]

  \node[main node] (1) at (.3,-.9) {};
  \node[main node] (2) at (0,0){};
  \node[main node] (3) at (.3,.9){};
  \node[main node] (4) at (1,0) {};
  \node[main node] (5) at (-.8,-.6) {};
  \node[main node] (6) at (-.8,.6) {};

  \path[every node/.style={font=\sffamily\small}]
   (1) edge node[above]{}(2)
   (2) edge node[above]{}(3)
   (2) edge node[above]{}(4)
   (2) edge node[above]{}(5)
   (2) edge node[above]{}(6);
\end{tikzpicture}
\hskip30mm
\begin{tikzpicture}[scale=1,
  thick,main node/.style={circle,draw,font=\sffamily\bfseries,minimum size=3mm}]

  \node[main node] (1) at (0,0) {};
  \node[main node] (2) at (1/2,5/6){};
  \node[main node] (3) at (0,10/6){};
  \node[main node] (4) at (2,0) {};
  \node[main node] (5) at (3/2,5/6) {};
  \node[main node] (6) at (2,10/6) {};

  \path[every node/.style={font=\sffamily\small}]
   (1) edge node[above]{}(2)
   (2) edge node[above]{}(3)
   (2) edge node[above]{}(5)
   (4) edge node[above]{}(5)
   (5) edge node[above]{}(6);
\end{tikzpicture}

\end{center}

Here another candidate for a sufficient condition.

\begin{thm}{Question}\label{quest:all-tree}
Assume $F$ is a finite metric space that satisfies all tree comparisons.
Is it true that $F$ is isometric to a subset of an Alexandrov space with nonnegative curvature?
\end{thm}

Note that even for finite metric space the all tree comparison has to be checked for an infinite set of trees since one point of the space may be used as a label for several vertexes in the tree.

There is a chance that for 5-point and 6-point metric spaces, the condition in Question~\ref{quest:all-tree} is also necessary. 
However, since there are nonnegatively curved Riemannian manifolds that do not satisfy 4(1)-tree comparison, 
Theorem~\ref{thm:convexity} implies that this condition can not be necessary for 7-point metric spaces.

\medskip

For any metric space $X$ with an isometric group action $G\acts X$ with closed orbits the quotient map $X\to X/G$ is a submetry.
In particular, by Theorem~\ref{thm:hilbert-quotient}, if $G\acts \HH$ is an isometric action with closed orbits on the Hilbert space, then the quotient space $\HH/G$ satisfies all tree comparisons.

\begin{thm}{Question}
Assume $X$ is a metric space satisfying all tree comparisons.
Is it always possible to construct an isometric group action with closed orbits on the Hilbert space $G\acts \HH$ such that $X$ is isometric to a subset in $\HH/G$?
\end{thm}

\parbf{On matrix inequality.}
The comparison for monopolar trees has an algebraic corollary which was used Urs Lang and Viktor Schroeder in \cite{LS}, a similar inequality was used by Karl-Theodor  Sturm  in \cite{sturm}. 

Namely, given a point array $p,x_1,\dots,x_n$ in a metric space $X$ consider the $n{\times}n$-matrix $M$ with the components 
\[m_{i,j}=\tfrac12\cdot(|x_i-p|^2+|x_j-p|^2-|x_i-x_j|^2).\]
If the tree comparison for $p/x_1,\dots,x_n$ holds, then 
\[\bm{s}\cdot M\cdot \bm{s}^\top\ge 0\eqlbl{eq:sMs}\]
for any vector $\bm{s}=(s_1,\dots,s_n)$ with nonnegative components.

The converse does not hold; that is, for some point array $p,x_1,\dots,x_n$ in a metric space the inequality \ref{eq:sMs} might hold, while the tree comparison for $p/x_1x_2x_3x_4x_5$ does not.
(We do not know an explicit way to describe tree comparisons using a system of inequalities.)

\begin{wrapfigure}{r}{22 mm}
\vskip-0mm
\begin{tikzpicture}[scale=1,
  thick,main node/.style={circle,draw,font=\sffamily\bfseries,minimum size=3mm}, ]

    \node[main node] (0) at (0,0){};
     \node[main node] (1) at (.6,-.8) {};
     \node[main node] (2) at (.9,.3){};
     \node[main node] (3) at (0,1) {};
  \node[main node] (4) at (-.9,.3) {};
  \node[main node] (5) at (-.6,-.8) {};

  \path[every node/.style={font=\sffamily\small}]
  (1) edge node[above]{}(2)
  (2) edge node[above]{}(3)
  (3) edge node[above]{}(4)
  (4) edge node[above]{}(5)
  (5) edge node[above]{}(1)
  (1) edge[dashed] node[above]{}(3)
  (2) edge[dashed] node[above]{}(4)
  (3) edge[dashed] node[above]{}(5)
  (4) edge[dashed] node[above]{}(1)
  (5) edge[dashed] node[above]{}(2);
\end{tikzpicture}
\end{wrapfigure}

An example can be constructed by perturbing the configuration on the plane as on the diagram ---
if the diameter of diagram is 1, 
then increasing the distances between the pairs of points connected by dashed lines by $\eps=10^{-9}$ and decreasing  the distances between the pairs of points connected by sold lines by $\delta=10^{-6}$ does the job.
The obtained metric 6-point metric space satisfies the matrix inequality with center at each point, but does not satisfy the tree comparison with the pole at the central point.

Many necessary conditions on finite subsets of nonnegatively curved Alexandrov spaces are known;
in addition to the comparisons discussed above, 
let us mention the authors results in \cite{lebedeva-petrunin} and \cite{petrunin}
and the Markov type inequality proved by Shin-ichi Ohta in \cite{ohta};
see also the survey by Assaf Naor \cite{naor} and the references there in.

\parbf{On tree comparisons in length spaces.}
Note that if a tree $T$ is not a path, then it contains a tripod as subtree.
Therefore $T$-tree comparison implies Alexandrov comparison, in particular any complete length-metric space satisfying $T$-tree comparison is a nonnegatively curved Alexandrov space.

It is straightforward to generalize Theorem~\ref{thm:3(1)+2(2)} to  Alexandrov spaces; that is, we have the following theorem.

\begin{thm}{Theorem}
A complete length space $L$ satisfies  3(1) or 2(2)-tree comparison if and only if $L$ is a nonnegatively curved Alexandrov space.
\end{thm}

We expect that Theorem~\ref{thm:convexity} (after appropriate reformulation) can be also generalized to Alexandrov spaces --- the only obstacle we see is the proof of Proposition~\ref{prop:CTIL}.
Such a generalization would characterize length-metric spaces satisfying most of 4(1)-tree comparison (as well as most of bipolar comparisons).

The Alexandrov spaces that satisfy 4(1)-tree comparison remind the quotients of Riemannian manifolds by isometric group actions. 
For example we expect that if a finite dimensional Alexandrov space without boundary $A$ satisfies 4(1)-tree comparison, then the tangent space at any point $p\in A$ is a product of a Euclidean space $\mathbb{E}^k$ and a cone $K$ over space $\Sigma$ of diameter at most $\tfrac\pi2$ (the space $\Sigma$ might be empty, in this case $K$ is a one-point space).
In particular it implies that the set of all metric singularities of $A$ is an extremal subset, see \cite{perelman-petrunin}.
(For big branchy trees, the properties of spaces with the tree comparison should remind the quotients of Hilbert space even more.)

\parbf{On MTW.}
Recall that the condition (\textit{\ref{thm:convexity:MTW}}) in Theorem~\ref{thm:convexity} is named MTW$^{\not\perp}$.

Recall that for Riemannian manifolds, 3(1) and 2(2)-tree comparisons are equivalent to the nonnegative sectional curvature and 4(1) (as well as all $m$($n$)-tree comparison if $\max\{m,n\}\ge 4$) are equivalent to CTIL+MTW$^{\not\perp}$.
The meaning of tree comparisons for 3(2) and 3(3) remains unclear.
By theorems \ref{thm:3(1)+2(2)} and \ref{thm:convexity}, it is between Alexandrov comparison and CTIL+MTW$^{\not\perp}$.
The main result in \cite{lebedeva} implies that that 3(3)-tree comparison is strictly stronger than Alexandrov comparison.

\begin{thm}{Question}
What are the relations between 3(2) and 3(3)-tree comparisons,  MTW$^{\not\perp}$ and MTW for Riemannian manifolds with or without CTIL condition?
\end{thm}

For example, as we mentioned, MTW$^{\not\perp}$ is stronger than MTW, but we fail to show that it is strictly stronger.
In other words, we do not have an example of a (CTIL) Riemannian manifold that is MTW, but not MTW$^{\not\perp}$.
Another example: it might happen that 3(3)-tree comparison is equivalent to MTW+CTIL which would provide a synthetic description of these conditions.

We also do not know whether globalization theorem holds, for MTW$^{\not\perp}$ (or equivalently for 4(1)-tree comparison); in other words \emph{is it true that local MTW$^{\not\perp}$ implies global MTW$^{\not\perp}$?}

The following two questions are well known for MTW;
partial answers are given in \cite{MTW+CTIL+} and \cite{loeper} correspondingly.
The question might be easier for MTW$^{\not\perp}$.

\begin{thm}{Question}
Is it true that MTW (or MTW$^{\not\perp}$) implies CTIL?
\end{thm}

Note that if the globalization holds for MTW$^{\not\perp}$, then by Theorem~\ref{thm:convexity}, MTW$^{\not\perp}$ implies CTIL.

\begin{thm}{Question}
Is it true that CTIL+MTW (or CTIL+MTW$^{\not\perp}$) on a compact Riemannian manifold implies TCP?
\end{thm}

By Theorem~\ref{thm:convexity}, CTIL+MTW$^{\not\perp}$ is equivalent to 4(1)-tree comparison. 
Therefore the MTW$^{\not\perp}$-version of the last question can be reformulated the following way:

\begin{thm}{Question} Is it true that 4(1)-tree comparison on a compact Riemannian manifold implies TCP?
\end{thm}

\parbf{Acknowledgements.} 
We would like to thank an anonymous referee and Alexander Lytchak for thoughtful and constructive comments to the preliminary version of this paper.